\documentclass[final]{style_jp}
\usepackage{color,calc}
\usepackage{colortbl}
\usepackage{epsfig}
\usepackage{pstricks,pst-node,pst-text,pst-3d}

\begin{document}

\begin{frontmatter}

\title{A different way to read off the time $-$\\
A new idea for a binary clock}
\author{J\"org Pretz}
\ead{jorg.pretz@cern.ch}

\address{ Physikalisches Institut, 
 Universit\"at Bonn
 Nu{\ss}allee 12,
 53115 BONN, Germany}

\begin{abstract}
A new idea for a binary clock is presented.
It displays the time using a triangular array of 
15 bits. It is shown that such a geometric, 
triangular arrangement is only possible because our
system of time divisions is based on a sexagesimal system
in which the number of minutes in 12 hours
equals the  factorial of a natural number (720=6!).
\end{abstract}
\end{frontmatter}
\section{Introduction}
\begin{figure}[!t]
\begin{center}
\includegraphics[width=0.7\textwidth]{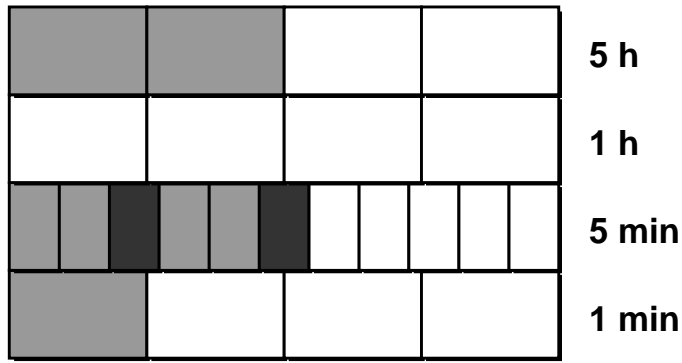}
\caption{Berlin clock.
The time displayed is 10:31.
Lamps lit correspond to dark areas.
\label{berlinuhr}}
\end{center}
\end{figure}

There are many ways to display the time.
For example the familiar analog display
with a dial and clock hands
or digital displays using numerals.
In addition there are binary displays which are a bit
more difficult to read.
A well known example is the so called Berlin clock
(also called ``set theory'' clock),
invented by  \mbox{Dieter} \mbox{Binninger} and first installed in 1975
on the Kurf\"urstendamm in Berlin. 
The display of this clock consists of bits realized by lamps
(or LEDs in a table version) which correspond to a certain
amount of time.
Lamps in the same row correspond to the same amount.
For the Berlin clock, every lamp in the top row
corresponds to 5 hours.
The corresponding values in the second row is 1 h.
In the third and fourth row the lamps represent
5 and 1 minute(s), respectively. 
For a better readability, every third lamp in the 
third row has a different color.
The time is then just given by adding the amounts of all 
lamps lit.
Fig.~\ref{berlinuhr} shows the Berlin clock.
The time displayed is 10:31.

\section{The new idea: A triangular clock}
\begin{figure}[!t]
\begin{center}
\includegraphics[width=0.6\textwidth]{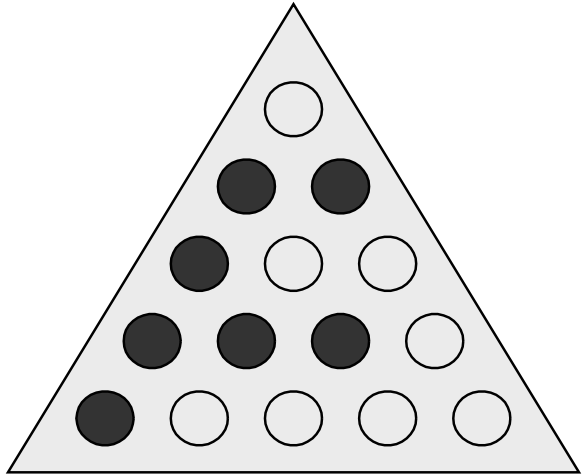}
\caption{Triangular Clock.
The time displayed is 4:49
(dark circle = lamp on).
\label{uhr}}
\end{center}
\end{figure}

A new, more esthetic and geometric, form of display will be presented now.
Here the lamps or LEDs are arranged in form of a triangle (see Fig.~\ref{uhr}).
Again, lamps in the same row correspond to the same amount of time.
The corresponding values are listed in the following table.
\begin{center}
\begin{tabular}{|c|r|}
\hline
row & value \\
\hline
1     &  6 hours \\
2     &  2 hours \\
3     &  30 minutes \\
4     &  6 minutes \\
5     &  1 minute \\
\hline
\end{tabular}
\end{center}

For the time displayed in Fig.~\ref{uhr} one finds:\\
\begin{center}
\begin{tabular}{rlcrr}
  $0 \times 6$&h   &=& 0h &\\
  $+2 \times 2$&h &=&  4h &\\
  $+1 \times 30$&min &=&  & 30 min\\
  $+3 \times 6$&min &=&   & 18 min\\
  $+1 \times 1$&min &=&   &  1 min\\
\hline
             &      & & 4h & 49 min \\          
\end{tabular}
\end{center}

\begin{figure}[!t]
\begin{center}
\includegraphics[width=0.32\textwidth]{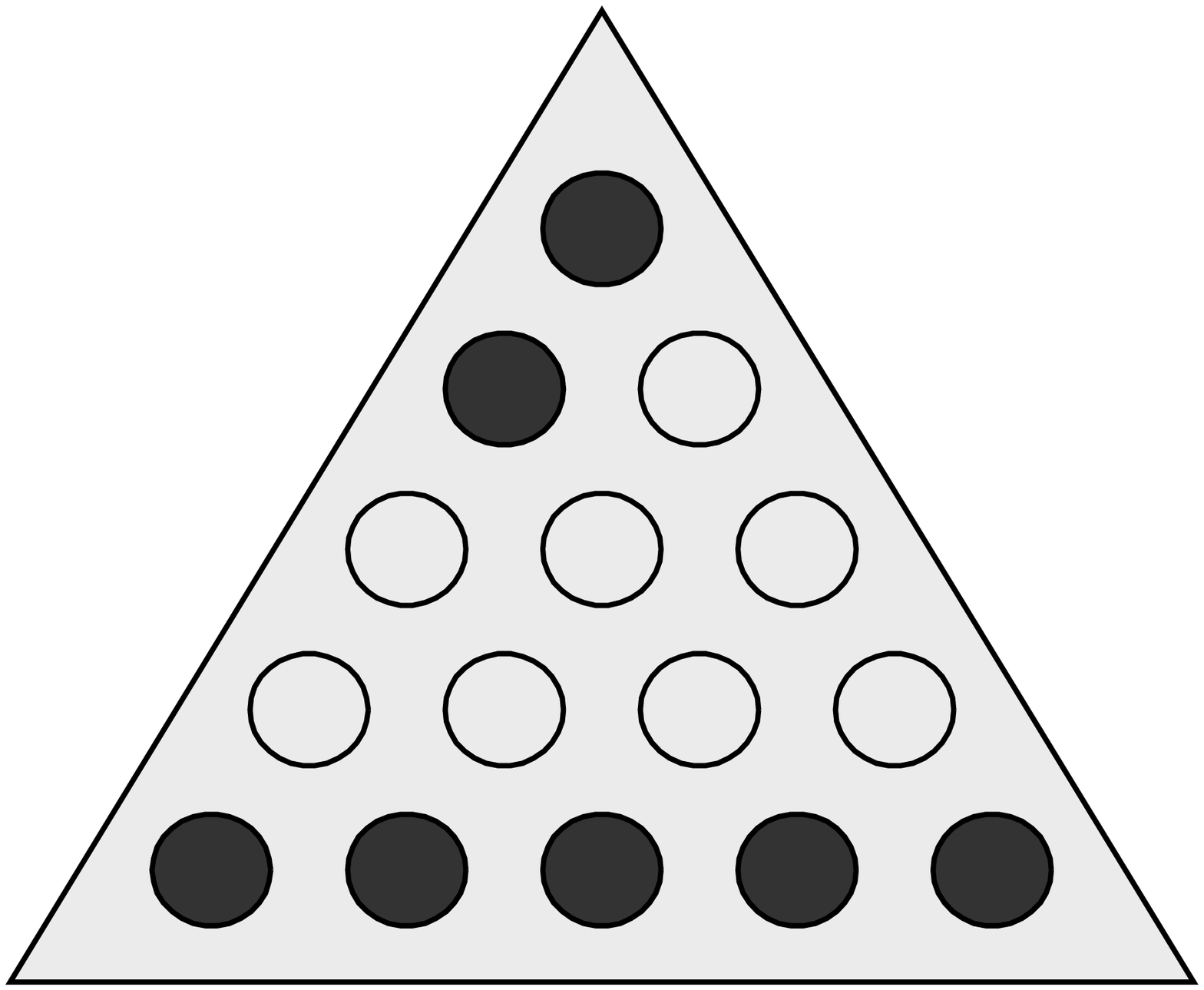} \vspace{2mm}
\includegraphics[width=0.32\textwidth]{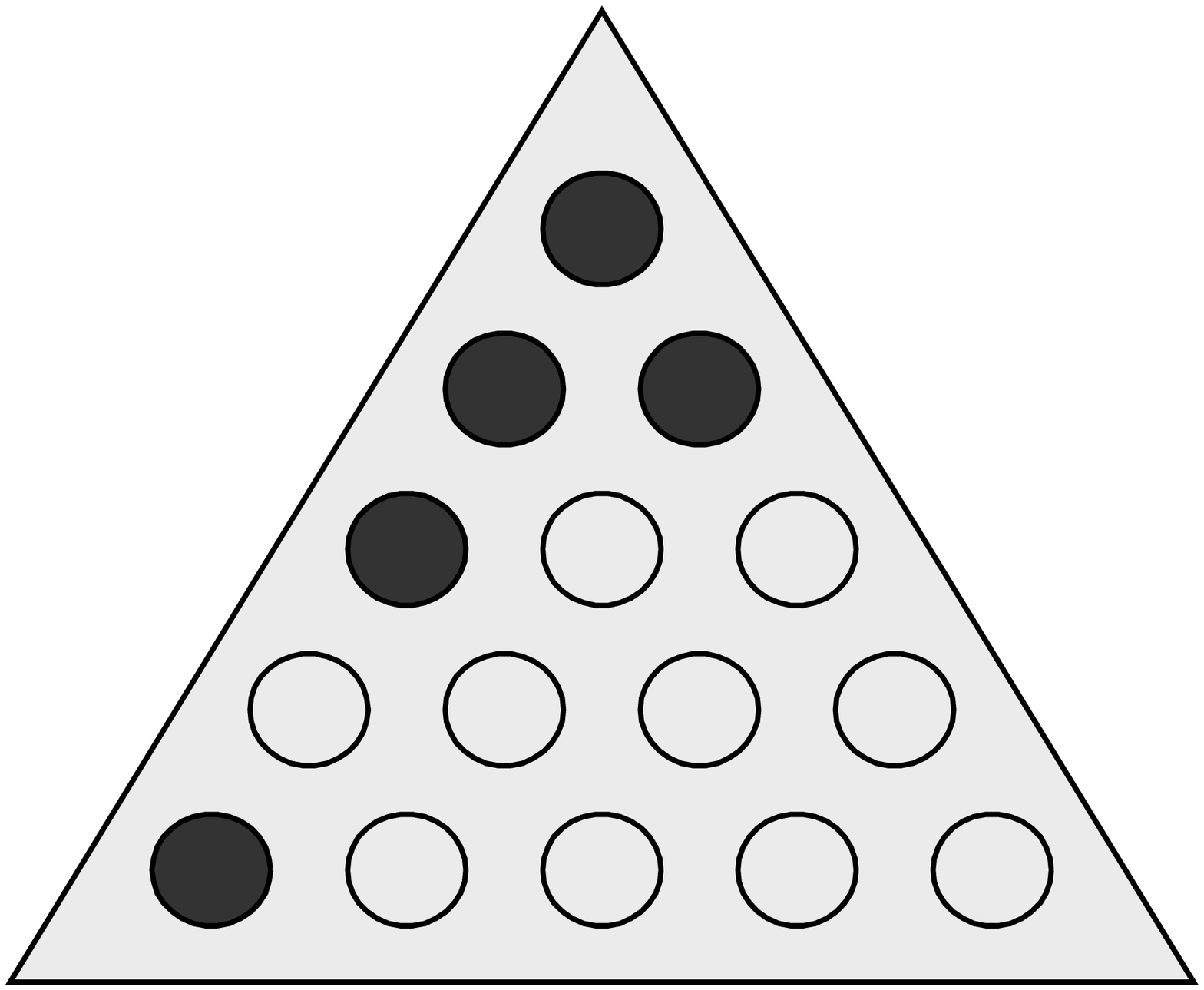}\vspace{2mm}
\includegraphics[width=0.32\textwidth]{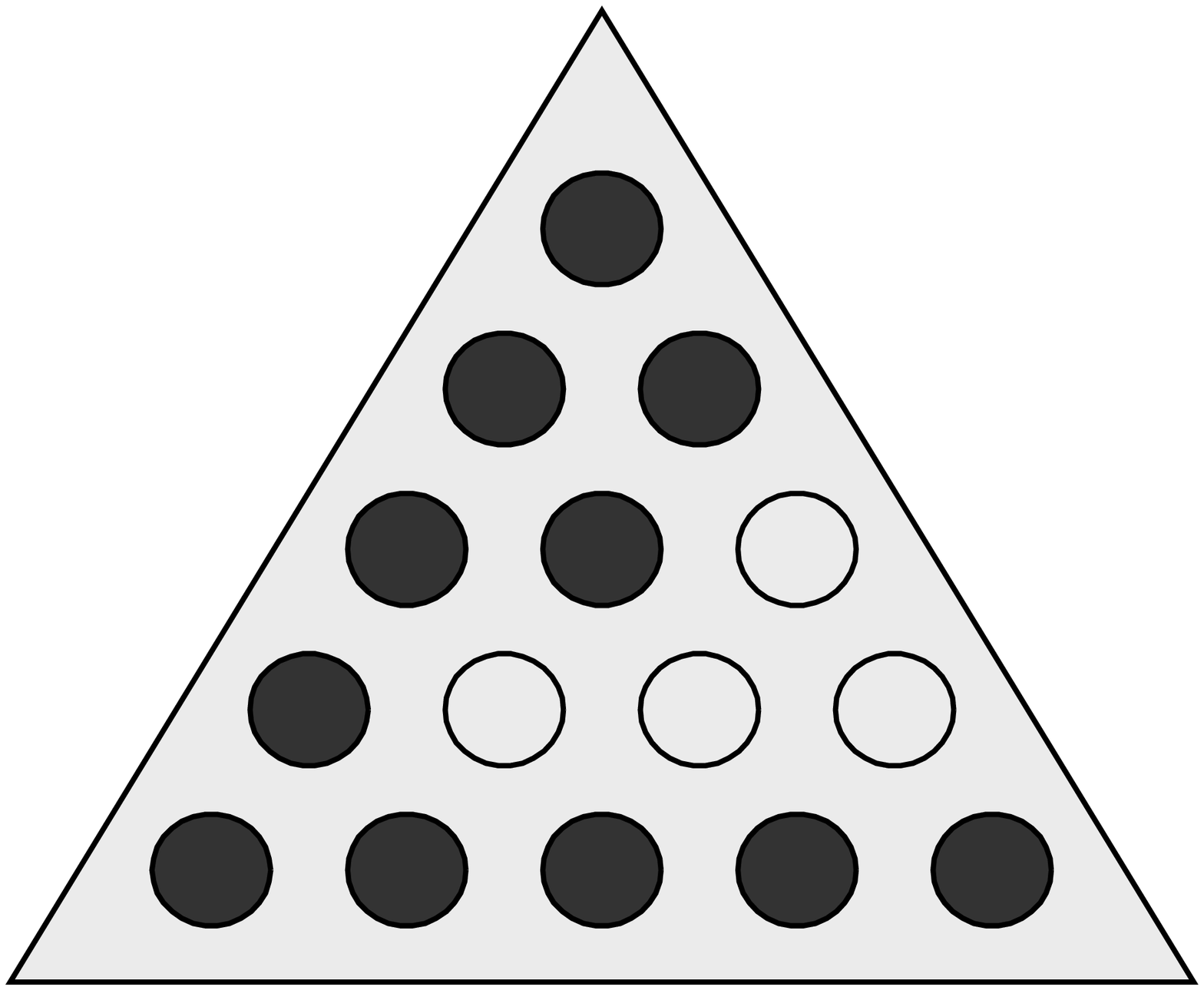}\vspace{2mm}
\caption{More examples. The times displayed are
08:05, 10:31 and 11:11.\label{moreex}}
\end{center}
\end{figure}

If all lamps are on, the time displayed is 
\[
6 \,\mbox{h} + 2 \times 2 \,\mbox{h} + 3 \times 30\, \mbox{min}  
+ 4 \times 6\, \mbox{min} + 5 \times 1 \,\mbox{min} = 11 \,\mbox{h} 59
\,\mbox{min} \, 
\]
which means that this triangular form with 5 rows of lamps
is perfectly suited for a 12 hour display.
Using two different colors, one can achieve a 24h display;
for example green for {\it am} and red for {\it pm}.
Fig.~\ref{moreex} shows more examples.
A Java applet of this clock can be found on
http://joerg.pretz.de. 

\section{Mathematical Background}
We now come to the question why our system of time divisions allows
such a triangular display?
First let us compare the Berlin and the triangular clock.
Common to both concepts is the fact that
the amount of time a certain lamp corresponds to in the $(n-1)$th row equals 
$(m_n+1)$ times the amount in the $n$th row
where $m_n$ is the number of lamps in the
$n$th row.
Here are two examples to illustrate this:
\begin{itemize}
   \item Berlin clock: 3rd row has $m_3 = 11$ lamps representing 5 minutes each.
    Thus for the 2nd row one finds $(11+1) \times 5$ min$ = 1$h.   
   \item triangular clock: 3rd row has $m_3 = 3$ lamps representing 30 minutes each.
    Thus for the 2nd row one finds $(3+1) \times 30$ min$ = 2$h.   
\end{itemize}
The main difference is that for the triangular clock one has
$m_n = n$, i.e.
in the $n$th row there are also $n$ lamps.
Only this allows the geometric, triangular arrangement of the lamps.

But what does $m_n=n$ mean?
A row with $n$ lamps allows to display $n+1$ states,
from all lamps off to all lamps on. 
(Note that in the way the display is used here, 
not all of the $2^n$ possible states are used because if a lamp is on, 
also the lamps left to it are on in the same row.)
For a triangular display with $n$ lamps
in the bottom row and the number of lamps decreasing by
one in every following row
the total number of states is thus
\[
(n+1) \times ((n-1)+1) \times ((n-2)+1) \times \dots \times (1+1) = (n+1)!  \, .
\]

Now, note that our system of time measurement
is based on divisions which have many factors,
e.g. $12=3 \times 4$ and $60 = 1 \times 2 \times 5 \times 6$. 
For a 12 hour display with a precision of one minute the number
of states one has to display is thus
\begin{eqnarray*}
{12 \times 60 \; \mbox{minutes}}  
&=&1 \times 2 \times 3 \times 4 \times 5 \times 6 \, \mbox{minutes} 
=6! \; \mbox{minutes} \, 
\end{eqnarray*}
which perfectly fits in a triangular display with five rows.
Thus the whole concept works because our system of time divisions is based 
on a sexagesimal system dating back to the Babylonian~\cite{neugebauer}
rather than a decimal system.

Fortunately an idea of the French Revolution to divide
the day in 10 hours of 100 minutes of 100 seconds didn't prevail~\cite{Carrigan}.
With such a division this triangular display would not be possible.


\begin{thebibliography}{99}
\bibitem{neugebauer} Otto Neugebauer, The Exact Sciences in Antiquity, Barnes \& Noble, New York, 1993
\bibitem{Carrigan} Richard A. Carrigan, Jr., Decimal Time, American Scientist, 66, 305 (78).
\end{thebibliography}
\end{document}